\documentclass[final,12pt]{colt2022_modified} 

\usepackage{latexsym}
\usepackage{calc}
\usepackage{amsfonts, listings, mathrsfs}

\def\nr{\par \noindent}

\def\sign{{\rm sign \,}}
\def\Def{\stackrel{\mathrm{def}}{=}}

\def\inter{{\rm int \,}}

\def\beq{\begin{equation}}
\def\eeq{\end{equation}}

\def\R{\mathbb{R}}

\newcommand{\qed}{\hfill $\Box$ \nr \medskip}

\def\argmin{\mathop{\rm argmin}}

\def\Argmax{\mathop{\rm Argmax}}

\def\ba{\begin{array}}
\def\ea{\end{array}}
\def\BT{\begin{theorem}}
\def\ET{\end{theorem}}
\def\BL{\begin{lemma}}
\def\EL{\end{lemma}}

\def\la{\langle}
\def\ra{\rangle}

\title[Lower bounds for minimizing regularized functions]
{Lower Complexity Bounds for Minimizing Regularized Functions}
\usepackage{times}

\coltauthor{%
 \Name{Nikita Doikov} \Email{nikita.doikov@uclouvain.be}\\
 \addr Université catholique de Louvain (UCLouvain), Louvain-la-Neuve, Belgium
	}

\begin{document}

\maketitle

\begin{abstract}%
In this paper, we establish lower bounds for the oracle
complexity of the first-order methods minimizing regularized convex functions. 
We consider the composite representation of the objective.
The smooth part has H\"older continuous gradient of degree $\nu \in [0, 1]$ and 
is accessible by a black-box local oracle. The composite part is a power of a norm.
We prove that the best possible rate
for the first-order methods in the large-scale setting
for Euclidean norms is of the order $O(k^{- p(1 + 3\nu) / (2(p - 1 - \nu))})$ 
for the functional residual,
where $k$ is the iteration counter and $p$ is the power of regularization.
Our formulation covers several cases, including computation
of the Cubically regularized Newton step by the first-order gradient methods, in which
case the rate becomes $O(k^{-6})$. It can be achieved by the Fast Gradient Method.
Thus, our result proves the latter rate to be optimal.
We also discover lower complexity bounds for non-Euclidean norms.
\end{abstract}


\section{Introduction}
\label{SectionIntro}

The modern complexity theory of Convex Optimization
originated in the book of~\cite{nemirovskii1983problem},
where the first lower bounds were constructed for many different 
classes of optimization problems and algorithms.
After this work, it became clear that
the best possible rate of convergence for a given family of methods 
is fundamentally limited by a particular problem class.

For example, in the case of minimizing the convex functions with bounded subgradient,
by black-box local methods, the best rate for decreasing the functional 
residual is $O(k^{-1/2})$,
assuming the problem dimension is sufficiently large.
This is the rate of the classical Subgradient Method (\cite{shor2012minimization}),
so the method is \textit{optimal}.
At the same time, for the differentiable convex 
functions with Lipschitz continuous gradient, the best rate is $O(k^{-2})$, 
which is much better. The optimal Fast Gradient Method
with this rate of convergence was developed by~\cite{nesterov1983method}.

During the last years, we started to see more and more examples of problem classes
that significantly vary the standard picture of complexity theory.
Thus, there were established
lower complexity bounds for the functions
with H\"older continuous gradient (w.r.t. different norms) by \cite{guzman2015lower},
for the functions with second- and high-order derivatives being Lipschitz continuous
by \cite{agarwal2018lower} and \cite{arjevani2019oracle}, 
for the relatively smooth functions by \cite{dragomir2021optimal}.
Different extensions of the complexity theory to the randomized algorithms
were presented by \cite{bubeck2019complexity},
\cite{diakonikolas2019lower}, and
\cite{garg2021near}.

Meanwhile, the framework of \textit{composite optimization},
that was proposed by \cite{nesterov2013gradient}, \cite{beck2009fast},
provided us with much more flexibility in formulation of problem parameters.
Indeed, in composite problems, we can add to the smooth part
of the objective additional terms such as indicator of a convex set,
or regularizers of a different kind. Consequently, smoothness characteristics of the objective 
(e.g. the Lipschitz constant) need not necessary be consistent with 
uniformly convex properties of the composite part. 

In this work, we investigate the complexity of
the optimization problems with composite regularizer 
that is given by a power of a norm.
While having a variety of examples for the problems of this type,
it appears that the lower complexity bounds are not really covered
in the literature.

Let us consider  convex minimization problems with the following structure,
\beq \label{MainProblem}
\ba{rcl}
\min\limits_{x \in \R^n} \Bigl\{  
F(x) & \Def & f(x) + \frac{\sigma}{p}\|x\|^p
\Bigr\},
\ea
\eeq
with some parameters $p \geq 2$, $\sigma > 0$, and
the standard Euclidean norm:
$$
\ba{rcl}
\|x\| & \Def & \biggl( 
\sum\limits_{i = 1}^n \bigl(  x^{(i)} \bigr)^2
\biggr)^{\! 1/2}.
\ea
$$
Function $f$ is the main source of complexity of solving~\eqref{MainProblem}.
We assume that its gradient is H\"older continuous,
for some $\nu \in [0, 1]$ and $H_{\nu} > 0$
\beq \label{HolderGrad}
\ba{rcl}
\| \nabla f(x) - \nabla f(y) \| & \leq & H_{\nu}\|x - y\|^{\nu},
\qquad \forall x, y \in \R^n.
\ea
\eeq

When $p = 2$, we obtain the classical $\ell_2$-regularization by the squared Euclidean norm.
Finding the minimum of~\eqref{MainProblem} is equivalent to computing the value
of the proximal operator of $f$ (\cite{moreau1965proximity}).
We know that the problem is \textit{strongly convex}.
The rate of convergence becomes \textit{linear}
for the functions with Lipschitz continuous gradient $(\nu = 1)$,
and $O(k^{-1})$ for bounded variation of gradient $(\nu = 0)$ (see, e.g. \cite{nesterov2018lectures}).
Thus, the cost of solving the regularized problem is much cheaper than minimizing $f$ solely.
It happens that these observations work for $p > 2$ as well.

For $p = 3$, minimization problem~\eqref{MainProblem} appears to be a subproblem
for computing one step of the Newton Method with Cubic Regularization~(\cite{nesterov2006cubic}).
The smooth part then is a convex second-order Taylor polynomial, which is a quadratic function,
$$
\ba{rcl}
f(x) & = & \frac{1}{2}\la Ax, x \ra - \la b, x \ra, \qquad A \succeq 0,
\ea
$$
and condition~\eqref{HolderGrad} is satisfied 
with $\nu = 1$ and constant $H_1 = \lambda_{\max}(A)$.
In the following years, the Cubic Newton Method received a substantial interest
with the development of inexact, stochastic, and adaptive schemes
(see the papers of \cite{cartis2011adaptive1,grapiglia2017regularized,doikov2018randomized,cartis2018global,hanzely2020stochastic}).
In all of these methods, it is required to solve the subproblem with 
our structure~\eqref{MainProblem}.
The Fast Gradient Methods with restarts that have the rate of convergence $O(k^{-6})$
for this problem were constructed by
\cite{roulet2017sharpness} and by \cite{nesterov2019inexact}
with applications to the Cubic Newton.
In our paper we justify that the latter rate is \textit{optimal}, i.e.
the best that can be achieved for the considered class of functions.
Note that in the methods of~\cite{grapiglia2017regularized,grapiglia2019accelerated},
it is also required to minimize the quadratic function
alongside the regularizer with \textit{arbitrary} $p \in [2, 3]$.

In a recent paper by~\cite{nesterov2022quartic},
there were considered new
second-order schemes based on regularization of \textit{degree} $p = 4$.
Then, the rate of the Fast Gradient Method for solving the corresponding subproblem is
$O(k^{-4})$. We prove that this rate is \textit{optimal}, by justifying the corresponding lower bound.

Finally, let us mention the framework of High-Order Proximal-Point Methods
(\cite{nesterov2021inexact})
that became influential for the development of \textit{super-fast} second-order schemes
(see also the works of \cite{kamzolov2020near,nesterov2021inexact2}).
By generalizing the construction of the classical Proximal-Point 
Algorithm (\cite{rockafellar1976monotone})
onto arbitrary order $p \geq 2$, the methods assume 
that the subproblem of the form~\eqref{MainProblem} can be solved
somehow at each iteration. Thus, it becomes an important question
to study the complexity of problems with arbitrary degree of regularization.

\smallskip

We are interested to analyse the worst-case behaviour
of the first-order algorithms on such problem classes. 
Each optimization algorithm $\mathcal{A}$ can be associated 
with a sequence of mappings
$$
\ba{rcl}
\mathcal{A} & = & (A_0, A_1, A_2, \dots).
\ea
$$
At each iteration $k \geq 0$, mapping $A_k$ takes as input 
\textit{oracle information} from the previous $k$ points, 
and returns the next point $x_{k+1}$. 
In other words, $x_1 := A_0(\varnothing)$, and
$$
\boxed{
\ba{rcl}
x_{k + 1} & := & A_k( \mathcal{I}(x_1), \dots, \mathcal{I}(x_{k}) ), \qquad k \geq 1.
\ea
}
$$
In this paper, we consider so called \textit{first-order local oracle} of $f$, that is
$$
\ba{rcl}
\mathcal{I}(x) & = & \{ f(x), \nabla f(x) \}.
\ea
$$
The parameters of the regularizer are directly available for the method.

When we provide a lower bound for the rate of convergence,
for a fixed method
we need to build an \textit{example of problem} from our class 
such that the residual after the first $k$ iterations is bounded from below:
$$
\ba{rcl}
F(x_k) - F^{*} & \geq & R_k,
\ea
$$
by some quantity $R_k$ that is called the \textit{risk}.
The inverse of it, i.e. the minimum number of iterations
required to have an $\varepsilon$-precision for the solution
is called the \textit{complexity}:
$$
\ba{rcl}
C_{\varepsilon} & = & \min \bigl\{ k \; : \; R_k < \varepsilon  \bigr\}.
\ea
$$

In the following table, we list the new lower bounds that we prove in our work
alongside the already known ones
for the risk and complexity.
The initial distance to the solution is denoted by~$D$.
Numerical constants are hidden.

\begin{center}
	{
		\small
		\centering
		\renewcommand{\arraystretch}{1.5}
		\begin{tabular}{ | l | c| c | c | }
			\multicolumn{4}{ c  }{  
				\begin{tabular}{c}
					{ Euclidean norms } 
					\\[7pt]
				\end{tabular} 
			}  
			\\
			\hline
			\textbf{ \scriptsize Specification } & 
			 \textbf{ \scriptsize Risk, }  $R_k$  & 
			 \textbf{ \scriptsize Complexity, }  $C_{\varepsilon}$  &
			 \textbf{ \scriptsize Reference} 
			\\
			\hline
			\begin{tabular}{@{}l@{}}
				\vspace*{10pt}
				$\;\, p > 1 + \nu$ \\[7pt]
			\end{tabular}
			 &
			$ \displaystyle \!\!\!\!\!
				\biggl( \, \frac{H_{\nu}}{\sigma^{\frac{1 + \nu}{p}} k^{\frac{1 + 3\nu}{2}} } 
				\, \biggr)^{\!\!\frac{p}{p - 1 - \nu}} \!\!\!\!
			$
			&
			$
			\displaystyle \!\!
			 \biggl( 
			 \frac{H_{\nu}}{\sigma^{\frac{1 + \nu}{p} } \varepsilon^{\frac{p - 1 - \nu}{p}} }
			 \biggr)^{\!\!\frac{2}{1 + 3\nu}} \!\!\!
			$
			&
			{ \scriptsize
			\textbf{Our result:}
			Theorem~\ref{TheoremMain} } \\
			\hline
			\begin{tabular}{@{}l@{}}
				\vspace*{10pt}
				\;{\scriptsize Lipschitz gradient}
				$(\nu = 1, p > 2)$ \\[7pt]
			\end{tabular}
			&
			$ \displaystyle \!\!\!\!\!
			\biggl( \, \frac{H_{1}}{\sigma^{2 / p} k^{2} } 
			\, \biggr)^{\! \frac{p}{p - 2}} \!\!\!\!
			$
			&
			$ \displaystyle
			\sqrt{ \frac{H_1}{\sigma^{2/p} \varepsilon^{(p - 2) / p}} }
			$
			&
			{\scriptsize $\triangle$} \\
			\hline
			
			\begin{tabular}{@{}l@{}}
				\vspace*{10pt}
				\;{\scriptsize Bounded gradient}
				$(\nu = 0, p \geq 2)$ \\[7pt]
			\end{tabular}
			&
			$ \displaystyle \!\!\!\!\!
			\biggl( \, \frac{H_{0}}{\sigma^{1 / p} k^{1/2} } 
			\, \biggr)^{\! \frac{p}{p - 1}} \!\!\!\!
			$
			&
			$ \displaystyle \!\!\!\!\!
			\biggl( \,  \frac{H_0}{\sigma^{1/p} \varepsilon^{(p - 1) / p}}   \, \biggr)^{\! 2}
			$
			&
			{ \scriptsize
			\cite{juditsky2014deterministic} } \\
			\hline
			\begin{tabular}{@{}l@{}}
				{ \scriptsize Strongly convex functions with }\\[-3pt]
				{ \scriptsize Lipschitz gradient } $(\nu = 1, p = 2)$ 
			\end{tabular} &
			$ \sigma D^2 \exp\Bigl( -k \sqrt{\frac{\sigma}{H_1}} \Bigr)$
			&
			$ \sqrt{\frac{H_1}{\sigma}} \log \Bigl( \frac{\sigma D^2}{\varepsilon} \Bigr)$
			&
			{ \scriptsize \cite{nesterov2018lectures} } \\
			\hline
			\begin{tabular}{@{}l@{}}
							\vspace*{5pt}
			 { \scriptsize No composite part } ($\sigma = 0$) \\[10pt]
			 \end{tabular}
			  &  
			 $ \displaystyle  \frac{H_{\nu}D^{1 + \nu}}{k^{(3\nu + 1) / 2}} $
			  &
			$\displaystyle  \Bigl( \frac{H_\nu D^{1 + \nu}}{\varepsilon} \Bigr)^{\! \frac{2}{3 \nu + 1} }$
			& \!\!\!\!\!
			{ \scriptsize \cite{guzman2015lower}  }
			\!\!\!\!\!
			\\
			\hline
	\end{tabular}
}
\end{center}

Note that the complexity lower bounds of the same order for similar problem classes
appeared for the first time in the work of \cite{nemirovskii1985optimal} (see their equation (1.21)).
However, it seems that the corresponding proof was never published for a wide audience.

Lower bounds for the cubic subproblem with a (possibly nonconvex) quadratic smooth part
were constructed by \cite{carmon2018analysis}. Comparing with that work,
we were interested in dependence of the complexity 
on the \textit{free} parameter $\sigma$ in our sublinear rates.
Thus, we admit arbitrary choice of the regularization constant.

Besides, one significant novelty of our analysis is the use of the \textit{composite} formulation
for~\eqref{MainProblem} that provides a complete flexibility for the regularizer
and smoothness parameters.
Exploring this possibility, we also construct lower bounds for the regularization
by non-Euclidean norms (see the table below).
The employing of arbitrary norms as regularizers for the methods with Taylor's polynomials of different
order was considered in a recent paper by \cite{gratton2021adaptive}.
One step of their second-order algorithm requires to have an inexact solution
to the problem of our form with $p = 3$.

\begin{center}
	{
		\small
		\centering
		\renewcommand{\arraystretch}{1.5}
		\begin{tabular}{ | l | c| c | c | }
			\multicolumn{4}{ c  }{  
				\begin{tabular}{c}
					{ Regularization by $\|\cdot\|_q$-norm, $q \geq 1$ } 
					\\[7pt]
				\end{tabular} 
			}  
			\\
			\hline
			\textbf{ \scriptsize Specification } & 
			\textbf{ \scriptsize Risk, }  $R_k$  & 
			\textbf{ \scriptsize Complexity, }  $C_{\varepsilon}$  &
			\textbf{ \scriptsize Reference} 
			\\
			\hline
			\begin{tabular}{@{}l@{}}
				\vspace*{10pt}
				$\;\, p > 1 + \nu$ \\[7pt]
			\end{tabular}
			&
			$ \displaystyle \!\!
			\biggl( \, \frac{H_{\nu}}{\sigma^{\frac{1 + \nu}{p}} k^{\frac{1 + \nu + \nu q}{q}} } 
			\, \biggr)^{\!\!\frac{p}{p - 1 - \nu}} \!\!\!\!
			$
			&
			$
			\displaystyle \!\!
			\biggl( 
			\frac{H_{\nu}}{\sigma^{\frac{1 + \nu}{p} } \varepsilon^{\frac{p - 1 - \nu}{p}} }
			\biggr)^{\!\!\frac{q}{1 + \nu + \nu q}} \!\!\!
			$
			&
			{ \scriptsize
				\textbf{Our result:}
				Theorem~\ref{TheoremGenMain} } \\
			\hline
			\begin{tabular}{@{}l@{}}
				\vspace*{10pt}
				\;{\scriptsize Lipschitz gradient}
				$(\nu = 1, p > 2)$ \\[7pt]
			\end{tabular}
			&
			$ \displaystyle \!\!\!\!\!
			\biggl( \, \frac{H_{1}}{\sigma^{2 / p} k^{ (2 + q) / q } } 
			\, \biggr)^{\! \frac{p}{p - 2}} \!\!\!\!
			$
			&
			$ \displaystyle
			\biggl( \, 
			\frac{H_1}{\sigma^{2 / p} \varepsilon^{(p - 2) / p}}
			\,\biggr)^{\! \frac{q}{2 + q}}
			$
			&
			{\scriptsize $\triangle$} \\
			\hline
			
			\begin{tabular}{@{}l@{}}
				\vspace*{10pt}
				\;{\scriptsize Bounded gradient}
				$(\nu = 0, p \geq 2)$ \\[7pt]
			\end{tabular}
			&
			$ \displaystyle \!\!\!\!\!
			\biggl( \, \frac{H_{0}}{\sigma^{1 / p} k^{1/q} } 
			\, \biggr)^{\! \frac{p}{p - 1}} \!\!\!\!
			$
			&
			$ \displaystyle \!\!\!\!\!
			\biggl( \, \frac{H_{0}}{\sigma^{1 / p} \varepsilon^{(p - 1) / p} } 
			\, \biggr)^{\! q} \!\!\!\!
			$
			&
			{\scriptsize $\triangle$} \\
			\hline
		\end{tabular}
	}
\end{center}

\bigskip

It appears that for $\ell_1$-norm we obtain the best lower bounds for the rate of convergence.
However, it is not clear how far these estimates from 
the upper bounds that can be achieved by some optimization schemes.
At the same time, the rates for the Euclidean case
are tight. They can be reached by the Fast Gradient Methods with restarts
(see, e.g. the monograph of \cite{d2021acceleration}).

The rest of the paper is organized as follows.
In Section~\ref{SectionSmoothing} we review the basic properties of the local smoothing.
We use them in our analysis.
Section~\ref{SectionBound} contains the proof of our lower bounds for the Euclidean norm,
and in Section~\ref{SectionNonEuclidean} we discuss generalizations to non-Euclidean norms.
Section~\ref{SectionDiscussion} contains some final remarks.

\section{Local Smoothing}
\label{SectionSmoothing}

In this section, we summarise well known facts about the local smoothing,
that is the core component of the worst-case example.

Let $g: \R^n \to \R$ be a convex 1-Lipschitz functions, not-necessary differentiable. 
Thus 
$$
\ba{rcl}
g(x) - g(y) & \leq & \|x - y\|, \qquad \forall x, y \in \R^n.
\ea
$$
For a fixed parameter $\mu > 0$, we denote the local smoothing of $g$ by
$$
\ba{rcl}
\mathcal{S}_{\mu}[g](x) & \Def &
\min\limits_{y \in \R^n} \Bigl\{  g(y) + \frac{\mu}{2}\|y - x\|^2  \Bigr\},
\qquad x \in \R^n.
\ea
$$
Let us review the basic properties of the function $G(x) := S_{\mu}[g](x)$.
The proofs can be found in the Appendix.

\BL \label{LemmaLipschitz}
Function $G$ is convex and differentiable. It is $1$-Lipschitz,
and its gradient is $\mu$-Lipschitz continuous. Thus, for all $x, y \in \R^n$:
\beq \label{LipFuncSmooth}
\ba{rcl}
G(x) - G(y) & \leq & \|x - y\|,
\ea
\eeq
and
\beq \label{LipGradSmooth}
\ba{rcl}
\| \nabla G(x) - \nabla G(y) \| & \leq & \mu \|x - y \|.
\ea
\eeq
\EL

As a direct consequence of this lemma, we can conclude that $G$ has 
a H\"older continuous gradient for \textit{any} $\nu \in [0, 1]$ with
constant 
\beq \label{H_nu_constant}
\ba{rcl}
H_{\nu}(G) & = & 2^{1 - \nu} \mu^{\nu}.
\ea
\eeq
Indeed, \eqref{LipFuncSmooth} implies that the gradients are bounded: 
$\| \nabla G(x) \| \leq 1$, $\forall x \in \R^n$. Hence,
$$
\ba{rcl}
\| \nabla G(x) - \nabla G(y) \|
& = & 
\| \nabla G(x) - \nabla G(y) \|^{\nu} \cdot \| \nabla G(x) - \nabla G(y) \|^{1 - \nu} \\
\\
& \leq &
\mu^{\nu} \|x - y\|^{\nu} \cdot 2^{1 - \nu}
\;\; = \;\; H_{\nu}(G)\|x - y\|^{\nu}, \qquad \forall x, y \in \R^n.
\ea
$$

The following lemma shows that $G$ is also close to $g$.

\BL \label{LemmaClose}
For any $x \in \R^n$, it holds
\beq \label{G_close_to_g}
\ba{rcl}
g(x) & \geq & G(x) \;\; \geq \;\; g(x) - \frac{1}{2\mu}.
\ea
\eeq
\EL

\begin{figure}[h]
	\centering
	\includegraphics[width=0.50\textwidth ]{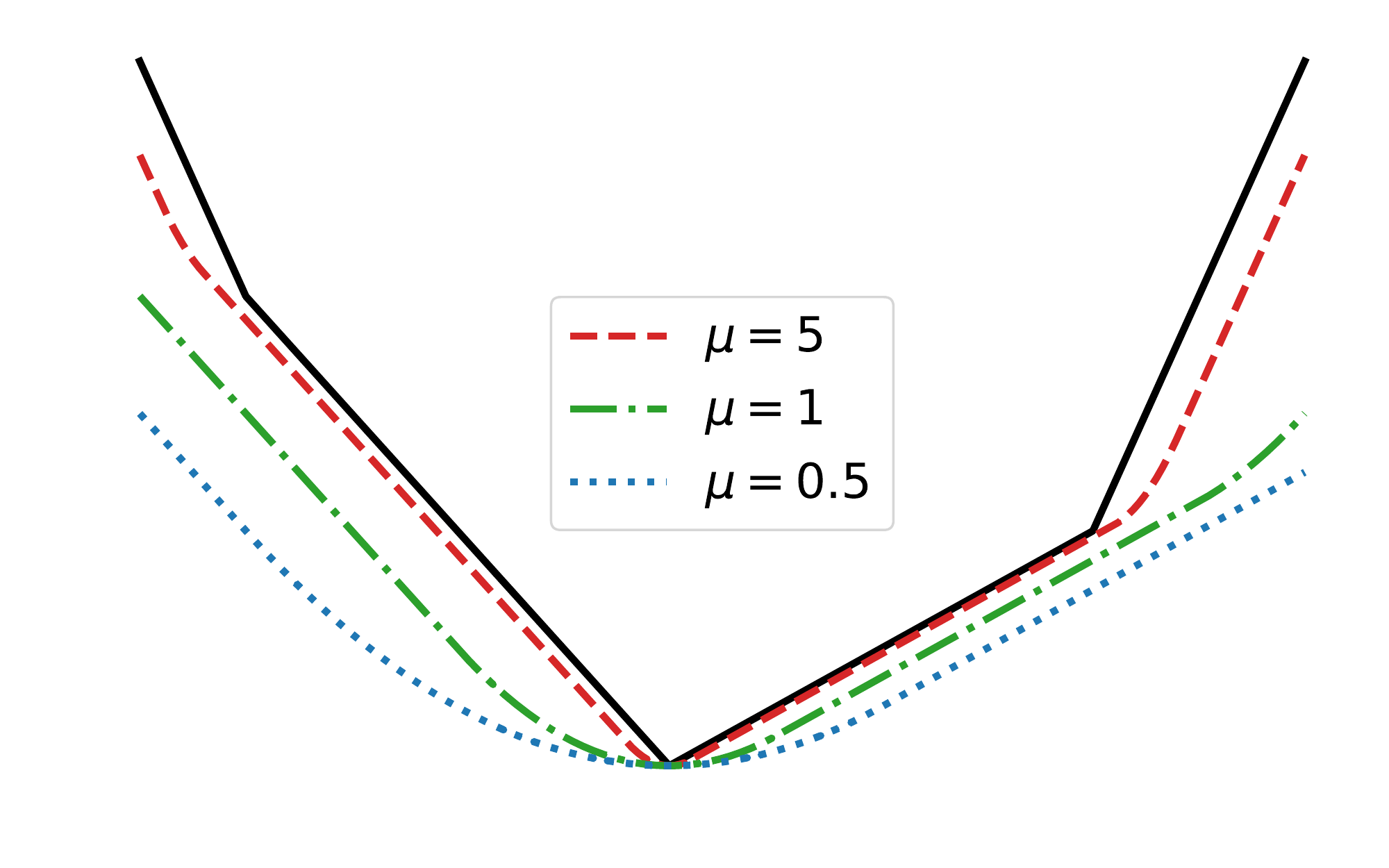}
	\caption{ \small Local smoothing of a piecewise linear function.}
	\label{Fig:Smoothing}
\end{figure}

We see that parameter $\mu$ provides a trade-off between the quality of approximation
and the level of smoothness for the new function.
The effect of varying this parameter is shown in Figure~\ref{Fig:Smoothing}.

Finally, we state the following important property that $G$ depends on $g$ 
in \textit{a local way}:
the value and the gradient of $G$ at $x$ depend only on the restriction
$g|_{U_{x, \mu}}$ of $g$ onto the ball
$$
\ba{rcl}
U_{x, \mu} & \Def & 
\Bigl\{ 
y \in \R^n \; : \; 
\|y - x\| \leq \frac{2}{\mu}
\Bigr\}.
\ea
$$
More formally, the following proposition holds.
\BL \label{LemmaLocal}
Let $z := \argmin\limits_{y \in \R^n} 
\Bigl\{ \; g(y) + \frac{\mu}{2}\|y - x\|^2 \; \Bigr\}$
for some fixed $x \in \R^n$ and $\mu > 0$.
Then
$$
\ba{rcl}
x - \frac{1}{\mu}\nabla G(x) 
& \equiv & z \;\; = \;\;
\argmin\limits_{y \in \inter U_{x, \mu}}
\Bigl\{ 
\;
g(y) + \frac{\mu}{2}\|y - x\|^2
\;
\Bigr\}.
\ea
$$
\EL

Therefore, in order to compute the oracle information
$\{ G(x), \nabla G(x) \}$ at some given point $x$,
it is enough to have an access only to the \textit{neighbourhood}
of function $g$ around this point.
This fact is crucial for proving the lower bound.

\section{Lower Complexity Bound}
\label{SectionBound}

Now, employing the construction introduced in \cite{guzman2015lower},
we are going to provide a resisting oracle strategy to establish the 
lower complexity bound for problem~\eqref{MainProblem}. 

Let us fix arbitrary 
$\nu \in [0, 1]$ such that $\nu < p - 1$, and some constant $H_{\nu} > 0$.
We prove the following result.

\BT \label{TheoremMain}
For any $T$-step algorithm with $T \leq n$, 
there exists a convex differentiable function $f(\cdot)$
whose gradient is H\"older continuous of degree $\nu$ with constant $H_{\nu}$, such that
for problem~\eqref{MainProblem}, we have

\beq \label{MainResult}
\ba{rcl}
F(x_T) - F^{*}
& \geq & 
\bigl(  \frac{p - 1}{p} \bigr)^{\frac{(p - 1)(1 + \nu)}{p - 1 - \nu}}
 \bigl( \frac{1}{2} \bigr)^{ \frac{(2p - 1)(1 + \nu)}{p - 1 - \nu} } \cdot
\biggl(
\frac{ H_{\nu} }{  \sigma^{\frac{1 + \nu}{p} } T^{ \frac{1 + 3\nu}{2} }}
\biggr)^{\frac{p}{p - 1 - \nu}}.
\ea
\eeq
The size of the solution to that problem is bounded as
\beq \label{SolBound}
\ba{rcl}
\|x^{*}\| & \leq & 
\bigl(  3(p - 1) 2^{p - 3}  \bigr)^{\frac{1}{p}}
\cdot 
\Bigl(\frac{1}{2} \bigl( \frac{p - 1}{4p} \bigr)^{\nu} \Bigr)^{\frac{1}{p - 1 - \nu}}
\biggl( 
\frac{H_{\nu}}{\sigma T^{\frac{1 + 3\nu}{2}}}
\biggr)^{\frac{1}{p - 1 - \nu}}.
\ea
\eeq
\ET
\proof

\noindent
Let us fix some positive $\delta$,  $\mu$, and $\beta$.

Resisting oracle chooses a set of numbers $\xi_1, \dots, \xi_T \in \{-1, 1\}$
and a permutation $k \mapsto \alpha(k) \in \{1, 2, \dots, T\}$.

Having these parameters, we consider the family of convex functions, $x \in \R^n$:
$$
\ba{rcl}
g_t(x) & = & 
\max\limits_{1 \leq k \leq t} \Bigl[ \; \xi_k \la e_{\alpha(k) }, x\ra \, - \, (k - 1) \delta  \; \Bigr],
\qquad 1 \leq t \leq T,
\ea
$$
where $e_i$ is the standard basis vector in $\R^n$ (recall that we assume that $n \geq T$).
It is clear that all $g_t(\cdot)$ are Lipschitz continuous with constant $1$.

Then, we are going to take as our final objective the function
$$
\boxed{
	\ba{rcl}
	f(x) & := & \beta \mathcal{S}_{\mu}[ g_T ] (x)
	\ea
}
$$

Let us estimate the minimum $F^{*}$ and the size of the solution $\|x^{*}\|$.
Note that we have, for any $x \in \R^n$
\beq \label{FUpperBound}
\ba{rcl}
F(x) & = & 
\beta \mathcal{S}_{\mu}[g_T](x) + \frac{\sigma}{p}\|x\|^p
\;\; \overset{\eqref{G_close_to_g}}{\leq} \;\;
\beta g_T(x) + \frac{\sigma}{p}\|x\|^p \\
\\
& \leq & 
\mathcal{H}(x)
\;\; \Def \;\;
\max\limits_{1 \leq k \leq T} \beta \xi_k \la e_{\alpha(k)}, x \ra
+ \frac{\sigma}{p}\|x\|^p.
\ea
\eeq
At the same time, we can also bound $F$ by $\mathcal{H}$ from below,
as follows
\beq \label{FLowerBound}
\ba{rcl}
F(x) & \overset{\eqref{G_close_to_g}}{\geq} &
\beta g_T(x) + \frac{\sigma}{p}\|x\|^p - \frac{\beta}{2\mu} 
\;\; \geq \;\;
\mathcal{H}(x) - \frac{\beta}{2\mu} -  \beta (T - 1) \delta.
\ea
\eeq
We can compute the minimum of $\mathcal{H}(\cdot)$
by using the symmetry within the problem,
\beq \label{MinUpperBound}
\ba{rcl}
F^{*} & \overset{\eqref{FUpperBound}}{\leq} &
\mathcal{H}^* \;\; = \;\; \min\limits_{x \in \R^n} \mathcal{H}(x)
\;\; = \;\;
\min\limits_{\gamma > 0} \Bigl\{  -\beta \gamma 
+ \frac{\sigma}{p} T^{\frac{p}{2}} \gamma^{p}  \Bigr\} \\
\\
& = &
- \frac{p - 1}{p} \cdot \Bigl( 
\frac{\beta^p}{\sigma T^{p / 2}}
\Bigr)^{\frac{1}{p - 1}}.
\ea
\eeq
From the uniform convexity of the regularizer (see, e.g. Lemma 2.5 in \cite{doikov2021minimizing})
we can bound the size of the solution,
\beq \label{SolSize}
\ba{rcl}
\bigl( \frac{1}{2} \bigr)^{p - 2}  \frac{\sigma}{p}\|x^*\|^p
& \leq & 
F(0) - F^{*}
\;\; \overset{\eqref{FUpperBound}, \eqref{FLowerBound} }{\leq} \;\;
\mathcal{H}(0) - \mathcal{H}^{*}
+ \frac{\beta}{2\mu} +  \beta (T - 1) \delta \\
\\
& = & 
\frac{p - 1}{p} \cdot \Bigl( 
\frac{\beta^p}{\sigma T^{p / 2}}
\Bigr)^{\frac{1}{p - 1}}
+ \frac{\beta}{2\mu} +  \beta (T - 1) \delta.
\ea
\eeq

Now, let us present a strategy for choosing $\xi_k$ and $\alpha(k)$.

\begin{itemize}
	\item At first step, the algorithm returns point $x_1$,
	which does not depend on the objective.
	Let us pick
	$$
	\ba{rcl}
	\alpha(1) & \in & \Argmax\limits_{1 \leq k \leq T} |\la e_k,  x_1 \ra |.
	\ea
	$$
	In other words, $\alpha(1)$ is an index of 
	a maximal element (in absolute value) among first $T$ coordinates of $x_1$.
	Then, we specify $\xi_1 \in \{-1, 1\}$ in a way that
	$$
	\ba{rcl}
	\xi_1 \la e_{\alpha(1)}, x_1 \ra & = & 
	|\la e_{\alpha(1)}, x_1 \ra|,
	\ea
	$$
	hence $\xi_1 =  \sign(\la e_{\alpha(1)}, x_1 \ra)$.
	So,
	$$
	\ba{rcl}
	g_1(x) & = & \xi_1 \la e_{\alpha(1)}, x \ra.
	\ea
	$$
	
	\item At step $2 \leq t \leq T$, assume that we have built function $g_{t - 1}(x)$.
	Let $x_t$ be the point of the trajectory of $\mathcal{A}$ at iteration $t - 1$,
	applied to the current objective:
	$$
	\ba{rcl}
	x_t & := & A_{t - 1}( \mathcal{I}_{t - 1}(x_1), \dots, \mathcal{I}_{t - 1}(x_{t - 1})),
	\ea
	$$
	where $\mathcal{I}_{t - 1}(x) \Def \{ \beta \mathcal{S}_{\mu}[g_{t - 1}](x)  , 
	\beta \nabla \mathcal{S}_{\mu}[g_{t - 1}](x)   \}$.
	Let us choose as $\alpha(t)$ the index of a maximal element of $x_t$ (in absolute value) among
	first $T$ coordinates, except $\alpha(1), \dots, \alpha(t - 1)$. Thus,
	$$
	\ba{rcl}
	\alpha(t) & \in & 
	\Argmax\limits_{ \stackrel{1 \leq k \leq T}{k \notin \{\alpha(1), \dots, \alpha(t - 1) \}}} 
	|\la e_k, x_t \ra|.
	\ea
	$$
	Then, we specify $\xi_t \in \{-1, 1\}$ such that
	$$
	\ba{rcl}
	\xi_t \la e_{\alpha(t)}, x_t \ra & = & | \la e_{\alpha(t)}, x_t \ra |,
	\ea
	$$
	hence $\xi_t = \sign(\la e_{\alpha(t)}, x_t \ra)$, and thus we obtain the next $g_t(x)$.
\end{itemize}

We need to prove that for any $2 \leq t \leq T$, function $g_t(\cdot)$ coincides with $g_s(\cdot)$,
$1 \leq s < t$,
in the ball of radius $\frac{2}{\mu}$ around $x_s$:
\beq \label{LocalEquiv}
\ba{rcl}
g_t(x) & = & g_s(x), \qquad x \in U_{x_s, \mu}
\;\; \Def \;\; 
\Bigl\{ x \in \R^n \; : \; \|x - x_s\| \leq \frac{2}{\mu}  \Bigr\},
\ea
\eeq
which means that $g_t(\cdot)$ is \textit{indistinguishable} from $g_s(\cdot)$ during the first $s$ steps,
and
hence by Lemma~\ref{LemmaLocal} their smoothings
$\mathcal{S}_{\mu}[g_t]$ and $\mathcal{S}_{\mu}[g_s]$ are also indistinguishable for the algorithm.

Indeed,
$$
\ba{rcl}
g_t(x) & = & 
\max\Bigl\{  g_s(x), \; 
\max\limits_{s < k \leq t}\Bigl[\;  \xi_k\la e_{\alpha(k)}, x \ra \, - \, (k - 1) \delta    \;\Bigr] 
\Bigr\}.
\ea
$$
By the definition of $g_s(x)$, we have
$$
\ba{rcl}
g_s(x_s) & \geq & \Bigl[ \xi_k \la e_{\alpha(k)}, x_s \ra - (k - 1)\delta \Bigr] + \delta, 
\qquad s < k \leq t.
\ea
$$
Hence, due to the Lipschitz continuity, it holds for all $x$ such that 
$\|x - x_s\| \leq \frac{\delta}{2}$:
$$
\ba{rcl}
g_s(x) & \geq & \Bigl[ \xi_k \la e_{\alpha(k)}, x \ra - (k - 1) \delta \Bigr],
\qquad s < k \leq t.
\ea
$$
Choosing
\beq \label{DeltaChoice}
\boxed{
	\ba{rcl}
	\delta := \frac{4}{\mu}
	\ea
}
\eeq
we conclude that~\eqref{LocalEquiv} is true.

Thus, we established correctness for the construction of the resisting oracle. Namely, it holds,
for all $s \leq t \leq T$:
$$
\ba{rcl}
 \mathcal{I}_{t}(x_s) 
    & = &
 \mathcal{I}_T(x_s)
 \;\; \equiv \;\;
 \mathcal{I}(x_s),
\ea
$$
so the oracles $\mathcal{I}_{t}(\cdot)$ and $\mathcal{I}(\cdot)$
are identical along the trajectory of the method.

It remains to bound from below
the residual in the function value for the last iteration of the method.
From simple observations, we notice that
\beq \label{GTLastIter}
\ba{rcl}
g_T(x_T) & \geq & | x_T^{(\alpha(T))} | - (T - 1) \delta
\;\; \geq \;\; - (T - 1) \delta 
\;\; \overset{\eqref{DeltaChoice}}{=} \;\; -(T - 1) \frac{4}{\mu}.
\ea
\eeq
Hence, for the last function value, we have
\beq \label{FTLastIter}
\ba{rcl}
F(x_T) & \geq & f(x_T)
\;\; = \;\;
\beta \mathcal{S}_{\mu}[g_T](x_T)
\;\; \overset{ \eqref{G_close_to_g}, \eqref{GTLastIter}}{\geq} \;\; 
-(T - 1) \frac{4\beta}{\mu} - \frac{\beta}{2 \mu}
\;\; \geq \;\;
-\frac{4\beta T }{\mu}.
\ea
\eeq
Therefore, the residual can be bounded as follows,
\beq \label{ResBound}
\ba{rcl}
F(x_T) - F^{*}
& \overset{\eqref{FTLastIter}, \eqref{MinUpperBound}}{\geq} &
\frac{p - 1}{p} \cdot \Bigl( 
\frac{\beta^p}{\sigma T^{p / 2}}
\Bigr)^{\frac{1}{p - 1}}  - \frac{4\beta T }{\mu} 
\;\; = \;\; 
\frac{p - 1}{2p} \cdot \Bigl( 
\frac{\beta^p}{\sigma T^{p / 2}}
\Bigr)^{\frac{1}{p - 1}},
\ea
\eeq
where the last equality holds by the following choice of the smoothing parameter:
\beq \label{MuChoice}
\boxed{
\ba{rcl}
\mu & := & 
\frac{8p}{p - 1} \Bigl( 
\frac{\sigma T^{ (3p - 2) / 2 }}{\beta}
\Bigr)^{\frac{1}{p - 1}}
\ea
}
\eeq
The only parameter which remains to determine is $\beta$. 
We know that $f$ must have H\"older continuous gradient with constant
$$
\ba{rcl}
H_{\nu} & \overset{\eqref{H_nu_constant}}{=} &
\beta 2^{1 - \nu} \mu^{\nu}
\;\; \overset{\eqref{MuChoice}}{=} \;\;
2 \Bigl(  \frac{4p}{p - 1}  \Bigr)^{\nu}
\Bigl(  \sigma T^{(3p - 2) / 2}  \Bigr)^{\frac{\nu}{p - 1}}
\beta^{\frac{p - 1 - \nu}{p - 1}}.
\ea
$$
Therefore, we get the following value for the last free parameter:
\beq \label{BetaChoice}
\boxed{
\ba{rcl}
\beta & := & 
\Bigl( 
\frac{1}{2} \bigl(  \frac{p - 1}{4p} \bigr)^{\nu} H_{\nu}
\Bigr)^{\frac{p - 1}{p - 1 - \nu}}
\Bigl(  \sigma T^{(3p - 2) / 2} \Bigr)^{- \frac{\nu}{p - 1 - \nu}}
\ea
}
\eeq
Substituting this value into~\eqref{ResBound},
we obtain the final lower bound on the convergence rate~\eqref{MainResult}.

The bound~\eqref{SolBound} for the size of the solution follows immediately
from
$$
\ba{rcl}
\|x^{*}\| & \overset{\eqref{SolSize}}{\leq} & 
\biggl[ \frac{p 2^{p - 2}}{\sigma}
 \biggl( 
 \frac{p - 1}{p} \cdot \Bigl(  \frac{\beta^p}{\sigma T^{p / 2}} \Bigr)^{\frac{1}{p - 1}}
 + \frac{\beta}{2\mu} + \beta(T - 1) \delta
 \biggr)
\biggr]^{1 / p} \\
\\
& \overset{\eqref{DeltaChoice}}{\leq} &
\biggl[ \frac{p 2^{p - 2}}{\sigma}
\biggl( 
\frac{p - 1}{p} \cdot \Bigl(  \frac{\beta^p}{\sigma T^{p / 2}} \Bigr)^{\frac{1}{p - 1}}
+ \frac{4\beta T}{\mu}
\biggr)
\biggr]^{1 / p} \;\; \overset{\eqref{MuChoice}}{=} \;\;
\biggl[ \frac{3 (p - 1) 2^{p - 3}}{\sigma} 
\cdot \Bigl(  \frac{\beta^p}{\sigma T^{p / 2}} \Bigr)^{\frac{1}{p - 1}}
\biggr]^{1 / p},
\ea
$$
by substituting the expression for $\beta$ \eqref{BetaChoice}.
\qed

Let us consider the case $\nu = 1$. Then, Theorem \ref{TheoremMain}
states that for any $T$-step algorithm there is an instance of problem
with bounded solution,
\beq \label{ExplSolBound}
\ba{rcl}
\| x^{*} \| & \overset{\eqref{SolBound}}{\leq} & O \Bigl(  
\bigl( \frac{H_1}{\sigma T^2} \bigr)^p  \Bigr),
\ea
\eeq
such that the residual is bounded from below,
\beq \label{ExpLower}
\ba{rcl}
F(x_T) - F^{*} & \geq & 
O\Bigl( 
\bigl(  \frac{H_1}{\sigma^{2 / p} T^2}  \bigr)^{\frac{p}{p - 2}}
\Bigr).
\ea
\eeq
At the same time, the standard rate of the composite Fast Gradient Method~\cite{nesterov2013gradient}
starting from $x_1 = 0$ is
$$
\ba{rcl}
F(x_T) - F^{*} & \leq & O\Bigl( \frac{H_1 \|x^{*}\|^2}{T^2}  \Bigr)
\;\; \overset{\eqref{ExplSolBound}}{\leq} \;\;
O\Bigl( \bigl( \frac{H_1}{\sigma^{2/p} T^{2}} \bigr)^{\frac{p}{ p - 2}}  \Bigr),
\ea
$$
that matches the lower bound~\eqref{ExpLower}.

\newpage
\section{Regularization with Non-Euclidean Norms}
\label{SectionNonEuclidean}

In this section, we study a generalization 
of problem~\eqref{MainProblem} to non-Euclidean norms.
For arbitrary $q \geq 1$, we denote by $\| \cdot \|_q$
the standard $\ell_q$-norm:
$$
\ba{rcl}
\|x\|_q & \Def & \biggl(  
\sum\limits_{i = 1}^n 
\big| x^{(i)} \big|^{q}  \biggr)^{1/q}.
\ea
$$

Now, let us consider the problem:
\beq \label{GenMainProblem}
\ba{rcl}
\min\limits_{x \in \R^n} \Bigl\{ 
F(x) & \Def & f(x) + \frac{\sigma}{p}\|x\|^p_q
\Bigr\},
\ea
\eeq
for some $p \geq 2$, $\sigma > 0$, and convex differentiable $f$
that has H\"older continuous gradient
(w.r.t. the standard Euclidean norm\footnote{Hence, we use 
	different norms for the regularizer and for defining the smoothness properties of $f$. 
	Generalizations to non-Euclidean norms for $f$ are also possible, 
	see~(\cite{guzman2015lower}).}).
As before, we assume that $\nu < p - 1$.
We can prove the following result.

\BT \label{TheoremGenMain}
For any $T$-step algorithm with $T \leq n$, there exists
a convex differentiable function $f(\cdot)$ whose gradient
is H\"older continuous of degree $\nu$ with constant $H_{\nu}$,
such that for problem~\eqref{GenMainProblem}, we have
\beq \label{GenProblemLowerBound}
\ba{rcl}
F(x_T) - F^* & \geq & 
\bigl( \frac{p - 1}{p} \bigr)^{\frac{(p - 1)(1 + \nu)}{p - 1 - \nu}}
\bigl( \frac{1}{2} \bigr)^{\frac{(2p - 1)(1 + \nu)}{p - 1 - \nu}}
\cdot
\biggl(
\frac{ H_{\nu} }{  \sigma^{\frac{1 + \nu}{p} } T^{ \frac{1 + \nu + \nu q}{q} }}
\biggr)^{\frac{p}{p - 1 - \nu}}.
\ea
\eeq
\ET
\proof
The proof is very similar to that one of Theorem~\ref{TheoremMain}.
We use the same strategy for the resisting oracle, and the 
candidate for the smooth part is $f(x) := \beta \mathcal{S}_{\mu}[g_T](x)$,
with $\delta := \frac{4}{\mu}$ (see definition of $g_T(\cdot)$ in the previous theorem).

For the minimum, we have the following estimate,
$$
\ba{rcl}
F^{*} & = &  \min\limits_{x \in \R^n} F(x) 
\;\; \leq \;\;   \min\limits_{x \in \R^n} 
\biggl\{ \;  \max\limits_{1 \leq k \leq T} \beta \xi_k \la e_{\alpha(k)}, x \ra    
+ \frac{\sigma}{p}\|x\|_q^p
\; \biggr\} \\
\\
& = &
\min\limits_{\gamma > 0} 
\Bigl\{   -\beta \gamma + \frac{\sigma}{p} T^{\frac{p}{q}} \gamma^p  \Bigr\}
\;\; = \;\;
- \frac{p - 1}{p} \cdot \Bigl( 
\frac{\beta^p}{\sigma T^{p / q}}
\Bigr)^{\frac{1}{p - 1}}.
\ea
$$
Thus, we get the bound for the residual,
\beq \label{GenMainResBound}
\ba{rcl}
F(x_T) - F^{*} & \geq & 
\frac{p - 1}{p} \cdot  \Bigl( 
\frac{\beta^p}{\sigma T^{p / q}}
\Bigr)^{\frac{1}{p - 1}} - \frac{4\beta T}{\mu}
\;\; = \;\; 
\frac{p - 1}{2p} \cdot \Bigl( \frac{\beta^p}{\sigma T^{p / q}} \Bigr)^{ \frac{1}{p - 1}},
\ea
\eeq
where we made the following choice of the smoothing parameter to balance the two terms,
$$
\boxed{
\ba{rcl}
\mu & := & \frac{8p}{p - 1} \Bigl( 
\frac{\sigma T^{\frac{p + qp - q}{q}}}{\beta}
\Bigr)^{\frac{1}{p - 1}}
\ea
}
$$
Substituting this expression into the equation for the constant of H\"older continuity,
$H_{\nu} = \beta 2^{1 - \nu} \mu^{\nu}$, we get the value for the last unknown parameter,
$$
\boxed{
\ba{rcl}
\beta & := & 
\Bigl( 
\frac{1}{2} \bigl(  \frac{p - 1}{4p}  \bigr)^{\nu} H_{\nu}
\Bigr)^{\frac{p - 1}{p - 1 - \nu}}
\Bigl( \sigma T^{\frac{p + qp - q}{q}} \Bigr)^{-\frac{\nu}{p - 1 - \nu}}
\ea
}
$$
Plugging it into~\eqref{GenMainResBound} completes the proof.
\qed

\begin{figure}[h]
	\centering
	\includegraphics[width=0.28\textwidth ]{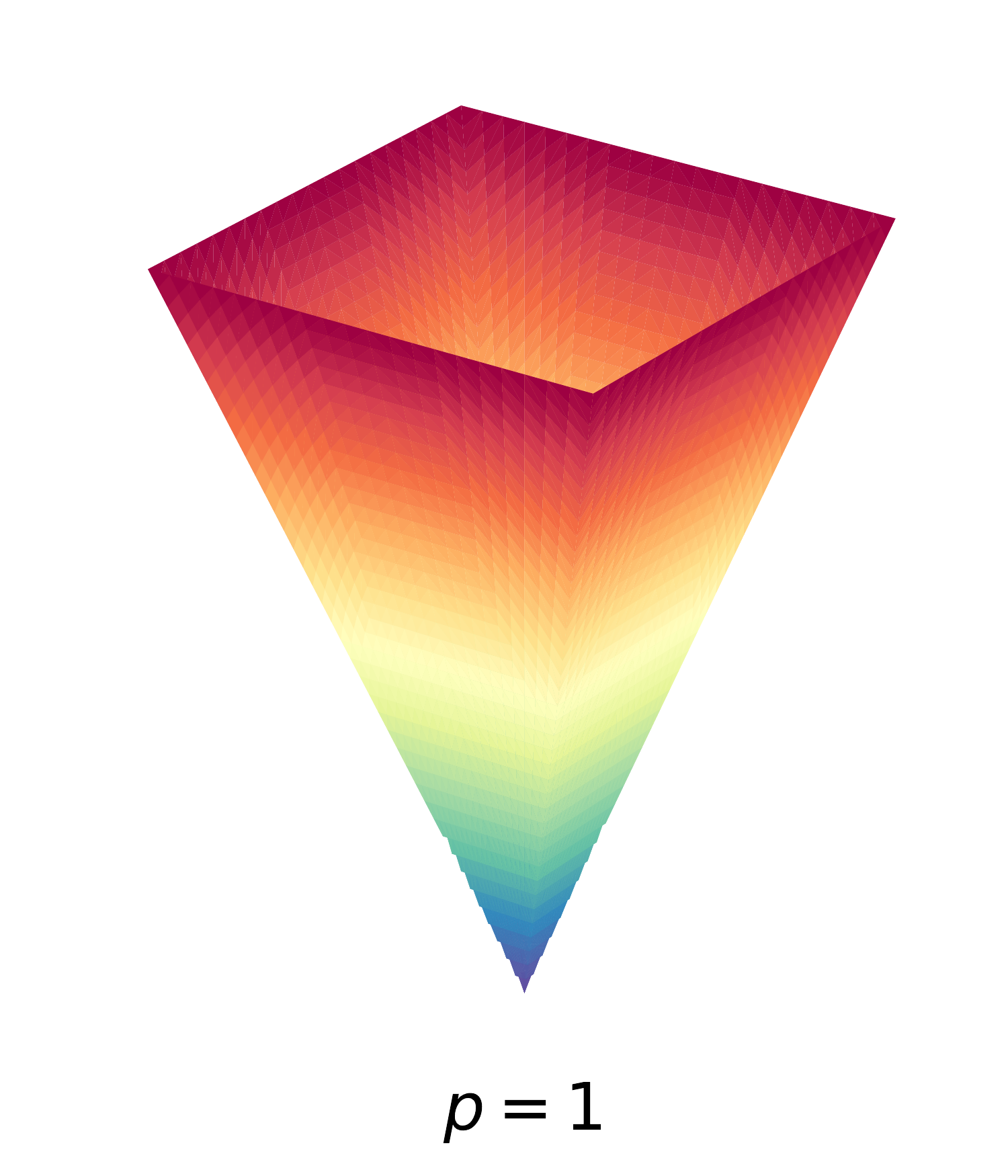}
	\includegraphics[width=0.28\textwidth ]{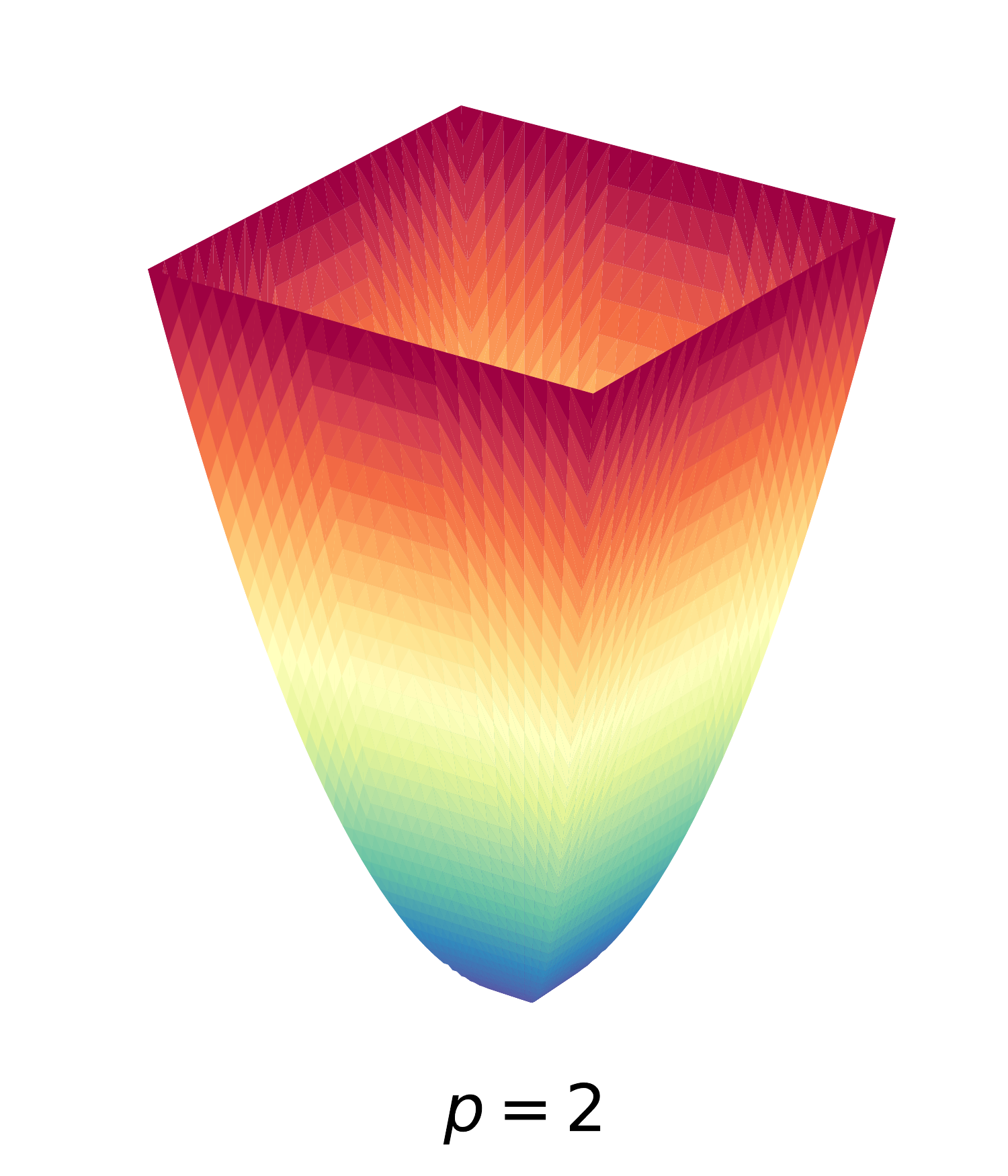}
	\includegraphics[width=0.28\textwidth ]{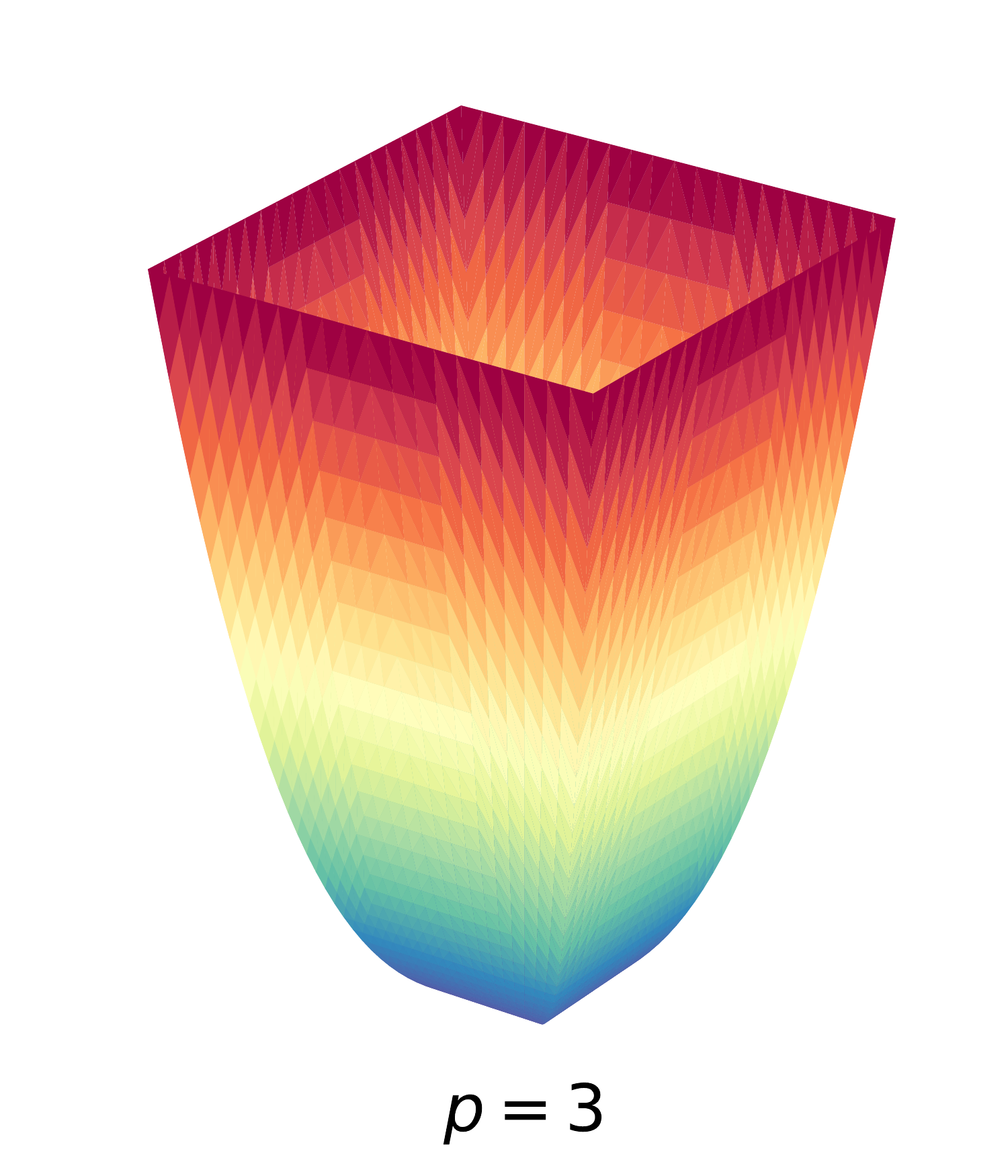}
	\caption{ \small Level sets of $\|\cdot\|_1^p$ for different powers.}
	\label{Fig:Norms}
\end{figure}

\section{Discussion}
\label{SectionDiscussion}

In this paper, we proved new lower bounds on the 
oracle complexity for minimizing regularized convex functions
by first-order methods.
As a particular case, we established that the 
\textit{best possible} rate of convergence for minimizing a function with Lipschitz continuous gradient
 regularized by \textit{cube} of the Euclidean norm is $O(k^{-6})$,
and for the \textit{fourth power} of the norm we have $O(k^{-4})$.
We know that the Fast Gradient Method achieves these rates.
It remains to be an interesting open question: whether we can construct the lower bounds
with a quadratic function for the smooth part.

Another interesting observation that we discovered is a change in the rate of convergence
for the regularizers with different norms. Thus, by using the third power of $\ell_1$-norm,
we obtain the lower bound of order $O(k^{-9})$. 
However, it is not clear how far this estimate from the upper bound
that can be reached by some optimization method.
We keep these questions for further investigation.
At the same time, taking a power of $\ell_1$-norm might preserve its sparsity properties (see Figure~\ref{Fig:Norms}),
which is desirable for applications.

\acks{This project has received funding from the European Research Council (ERC) 
	under the European Union’s Horizon 2020 research and innovation programme (grant agreement No. 788368).}

\bibliography{bibliography}

\begin{thebibliography}{35}
\providecommand{\natexlab}[1]{#1}
\providecommand{\url}[1]{\texttt{#1}}
\expandafter\ifx\csname urlstyle\endcsname\relax
  \providecommand{\doi}[1]{doi: #1}\else
  \providecommand{\doi}{doi: \begingroup \urlstyle{rm}\Url}\fi

\bibitem[Agarwal and Hazan(2018)]{agarwal2018lower}
Naman Agarwal and Elad Hazan.
\newblock Lower bounds for higher-order convex optimization.
\newblock In \emph{Conference On Learning Theory}, pages 774--792. PMLR, 2018.

\bibitem[Arjevani et~al.(2019)Arjevani, Shamir, and Shiff]{arjevani2019oracle}
Yossi Arjevani, Ohad Shamir, and Ron Shiff.
\newblock Oracle complexity of second-order methods for smooth convex
  optimization.
\newblock \emph{Mathematical Programming}, 178\penalty0 (1-2):\penalty0
  327--360, 2019.

\bibitem[Bauschke et~al.(2011)Bauschke, Combettes, et~al.]{bauschke2011convex}
Heinz~H Bauschke, Patrick~L Combettes, et~al.
\newblock \emph{Convex analysis and monotone operator theory in Hilbert
  spaces}, volume 408.
\newblock Springer, 2011.

\bibitem[Beck and Teboulle(2009)]{beck2009fast}
Amir Beck and Marc Teboulle.
\newblock A fast iterative shrinkage-thresholding algorithm for linear inverse
  problems.
\newblock \emph{SIAM journal on imaging sciences}, 2\penalty0 (1):\penalty0
  183--202, 2009.

\bibitem[Bubeck et~al.(2019)Bubeck, Jiang, Lee, Li, and
  Sidford]{bubeck2019complexity}
S{\'e}bastien Bubeck, Qijia Jiang, Yin-Tat Lee, Yuanzhi Li, and Aaron Sidford.
\newblock Complexity of highly parallel non-smooth convex optimization.
\newblock \emph{Advances in Neural Information Processing Systems}, 32, 2019.

\bibitem[Carmon and Duchi(2018)]{carmon2018analysis}
Yair Carmon and John~C Duchi.
\newblock Analysis of {K}rylov subspace solutions of regularized non-convex
  quadratic problems.
\newblock \emph{Advances in Neural Information Processing Systems}, 31, 2018.

\bibitem[Cartis and Scheinberg(2018)]{cartis2018global}
Coralia Cartis and Katya Scheinberg.
\newblock Global convergence rate analysis of unconstrained optimization
  methods based on probabilistic models.
\newblock \emph{Mathematical Programming}, 169\penalty0 (2):\penalty0 337--375,
  2018.

\bibitem[Cartis et~al.(2011)Cartis, Gould, and Toint]{cartis2011adaptive1}
Coralia Cartis, Nicholas~IM Gould, and Philippe~L Toint.
\newblock Adaptive cubic regularisation methods for unconstrained optimization.
  {P}art {I}: motivation, convergence and numerical results.
\newblock \emph{Mathematical Programming}, 127\penalty0 (2):\penalty0 245--295,
  2011.

\bibitem[d'Aspremont et~al.(2021)d'Aspremont, Scieur, and
  Taylor]{d2021acceleration}
Alexandre d'Aspremont, Damien Scieur, and Adrien Taylor.
\newblock Acceleration methods.
\newblock \emph{arXiv preprint arXiv:2101.09545}, 2021.

\bibitem[Diakonikolas and Guzm{\'a}n(2019)]{diakonikolas2019lower}
Jelena Diakonikolas and Crist{\'o}bal Guzm{\'a}n.
\newblock Lower bounds for parallel and randomized convex optimization.
\newblock In \emph{Conference on Learning Theory}, pages 1132--1157. PMLR,
  2019.

\bibitem[Doikov and Nesterov(2021)]{doikov2021minimizing}
Nikita Doikov and Yurii Nesterov.
\newblock Minimizing uniformly convex functions by cubic regularization of
  {N}ewton method.
\newblock \emph{Journal of Optimization Theory and Applications}, pages 1--23,
  2021.

\bibitem[Doikov and Richt{\'a}rik(2018)]{doikov2018randomized}
Nikita Doikov and Peter Richt{\'a}rik.
\newblock Randomized block cubic {N}ewton method.
\newblock In \emph{International Conference on Machine Learning}, pages
  1289--1297, 2018.

\bibitem[Dragomir et~al.(2021)Dragomir, Taylor, d’Aspremont, and
  Bolte]{dragomir2021optimal}
Radu-Alexandru Dragomir, Adrien~B Taylor, Alexandre d’Aspremont, and
  J{\'e}r{\^o}me Bolte.
\newblock Optimal complexity and certification of {B}regman first-order
  methods.
\newblock \emph{Mathematical Programming}, pages 1--43, 2021.

\bibitem[Garg et~al.(2021)Garg, Kothari, Netrapalli, and Sherif]{garg2021near}
Ankit Garg, Robin Kothari, Praneeth Netrapalli, and Suhail Sherif.
\newblock Near-optimal lower bounds for convex optimization for all orders of
  smoothness.
\newblock \emph{Advances in Neural Information Processing Systems}, 34, 2021.

\bibitem[Grapiglia and Nesterov(2017)]{grapiglia2017regularized}
Geovani~N Grapiglia and Yurii Nesterov.
\newblock Regularized {N}ewton methods for minimizing functions with
  {H}\"{o}lder continuous {H}essians.
\newblock \emph{SIAM Journal on Optimization}, 27\penalty0 (1):\penalty0
  478--506, 2017.

\bibitem[Grapiglia and Nesterov(2019)]{grapiglia2019accelerated}
Geovani~N Grapiglia and Yurii Nesterov.
\newblock Accelerated regularized {N}ewton methods for minimizing composite
  convex functions.
\newblock \emph{SIAM Journal on Optimization}, 29\penalty0 (1):\penalty0
  77--99, 2019.

\bibitem[Gratton and Toint(2021)]{gratton2021adaptive}
Serge Gratton and Philippe~L Toint.
\newblock Adaptive regularization minimization algorithms with non-smooth norms
  and euclidean curvature.
\newblock \emph{arXiv preprint arXiv:2105.07765}, 2021.

\bibitem[Guzm{\'a}n and Nemirovski(2015)]{guzman2015lower}
Crist{\'o}bal Guzm{\'a}n and Arkadi Nemirovski.
\newblock On lower complexity bounds for large-scale smooth convex
  optimization.
\newblock \emph{Journal of Complexity}, 31\penalty0 (1):\penalty0 1--14, 2015.

\bibitem[Hanzely et~al.(2020)Hanzely, Doikov, Richt{\'a}rik, and
  Nesterov]{hanzely2020stochastic}
Filip Hanzely, Nikita Doikov, Peter Richt{\'a}rik, and Yurii Nesterov.
\newblock Stochastic subspace cubic {N}ewton method.
\newblock In \emph{International Conference on Machine Learning}, pages
  4027--4038. PMLR, 2020.

\bibitem[Juditsky and Nesterov(2014)]{juditsky2014deterministic}
Anatoli Juditsky and Yurii Nesterov.
\newblock Deterministic and stochastic primal-dual subgradient algorithms for
  uniformly convex minimization.
\newblock \emph{Stochastic Systems}, 4\penalty0 (1):\penalty0 44--80, 2014.

\bibitem[Kamzolov and Gasnikov(2020)]{kamzolov2020near}
Dmitry Kamzolov and Alexander Gasnikov.
\newblock Near-optimal hyperfast second-order method for convex optimization
  and its sliding.
\newblock \emph{arXiv preprint arXiv:2002.09050}, 2020.

\bibitem[Moreau(1965)]{moreau1965proximity}
Jean-Jacques Moreau.
\newblock Proximity{\'e} and dualit{\'e} in a hilbertian space.
\newblock \emph{Bulletin of the Math{\'e}matic Society of France}, 93:\penalty0
  273--299, 1965.

\bibitem[Nemirovski and Nesterov(1985)]{nemirovskii1985optimal}
Arkadi Nemirovski and Yurii Nesterov.
\newblock Optimal methods of smooth convex minimization.
\newblock \emph{USSR Computational Mathematics and Mathematical Physics},
  25\penalty0 (2):\penalty0 21--30, 1985.

\bibitem[Nemirovski and Yudin(1983)]{nemirovskii1983problem}
Arkadi Nemirovski and David Yudin.
\newblock Problem complexity and method efficiency in optimization.
\newblock 1983.

\bibitem[Nesterov(1983)]{nesterov1983method}
Yurii Nesterov.
\newblock A method for solving the convex programming problem with convergence
  rate {O}(1/k\^{}2). [in {R}ussian].
\newblock In \emph{Dokl. akad. nauk Sssr}, volume 269, pages 543--547, 1983.

\bibitem[Nesterov(2013)]{nesterov2013gradient}
Yurii Nesterov.
\newblock Gradient methods for minimizing composite functions.
\newblock \emph{Mathematical Programming}, 140\penalty0 (1):\penalty0 125--161,
  2013.

\bibitem[Nesterov(2018)]{nesterov2018lectures}
Yurii Nesterov.
\newblock \emph{Lectures on convex optimization}, volume 137.
\newblock Springer, 2018.

\bibitem[Nesterov(2019)]{nesterov2019inexact}
Yurii Nesterov.
\newblock Inexact basic tensor methods.
\newblock \emph{CORE Discussion Papers}, 23:\penalty0 2019, 2019.

\bibitem[Nesterov(2021{\natexlab{a}})]{nesterov2021inexact}
Yurii Nesterov.
\newblock Inexact accelerated high-order proximal-point methods.
\newblock \emph{Mathematical Programming}, pages 1--26, 2021{\natexlab{a}}.

\bibitem[Nesterov(2021{\natexlab{b}})]{nesterov2021inexact2}
Yurii Nesterov.
\newblock Inexact high-order proximal-point methods with auxiliary search
  procedure.
\newblock \emph{SIAM Journal on Optimization}, 31\penalty0 (4):\penalty0
  2807--2828, 2021{\natexlab{b}}.

\bibitem[Nesterov(2022)]{nesterov2022quartic}
Yurii Nesterov.
\newblock Quartic regularity.
\newblock \emph{arXiv preprint arXiv:2201.04852}, 2022.

\bibitem[Nesterov and Polyak(2006)]{nesterov2006cubic}
Yurii Nesterov and Boris Polyak.
\newblock Cubic regularization of {N}ewton's method and its global performance.
\newblock \emph{Mathematical Programming}, 108\penalty0 (1):\penalty0 177--205,
  2006.

\bibitem[Rockafellar(1976)]{rockafellar1976monotone}
R~Tyrrell Rockafellar.
\newblock Monotone operators and the proximal point algorithm.
\newblock \emph{SIAM journal on control and optimization}, 14\penalty0
  (5):\penalty0 877--898, 1976.

\bibitem[Roulet and d'Aspremont(2017)]{roulet2017sharpness}
Vincent Roulet and Alexandre d'Aspremont.
\newblock Sharpness, restart and acceleration.
\newblock \emph{Advances in Neural Information Processing Systems}, 30, 2017.

\bibitem[Shor(2012)]{shor2012minimization}
Naum~Zuselevich Shor.
\newblock \emph{Minimization methods for non-differentiable functions},
  volume~3.
\newblock Springer Science \& Business Media, 2012.

\end{thebibliography}

\newpage
\appendix

\section{Properties of Local Smoothing}

In this section, we provide the proofs for basic properties
of local smoothing (Section~\ref{SectionSmoothing}).

For any fixed $x \in \R^n$, let us denote by $z_x$
the solution to the smoothing problem, i.e.
$$
\ba{rcl}
z_x & \Def & \argmin\limits_{y \in \R^n}
\Bigl\{ g(y) + \frac{\mu}{2}\|y - x\|^2 \Bigr\}.
\ea
$$
Hence, $z_x$ satisfies 
the following optimality condition (see, e.g. Theorem 3.1.23 in \cite{nesterov2018lectures}),
\beq \label{ApZStat}
\ba{rcl}
\la \mu(z_x - x), y - z_x \ra + g(y) & \geq & g(z_x), \qquad \forall y \in \R^n.
\ea
\eeq
In other words, $g'(z_x) \Def  -\mu(z_x - x)  \, \in \, \partial g(z_x)$.

\subsection{Proof of Lemma~\ref{LemmaLipschitz}}

We know a simple formula for the gradient of $G$, that is
$$
\boxed{
\ba{rcl}
\nabla G(x) & = & -\mu(z_x - x)
\ea
}
$$
(see, e.g. Proposition 12.30 in \cite{bauschke2011convex}).
Therefore, since $g$ is $1$-Lipschitz, we have
$$
\ba{rcl}
\| \nabla G(x) \| & = & \| g'(z_x) \| \;\; \leq \;\; 1, \qquad \forall x \in \R^n,
\ea
$$
and hence $G$ is also $1$-Lipschitz.

Now, let us fix two arbitrary points $x, y \in \R^n$. Then,
$$
\ba{rcl}
0 & \leq & \la g'(z_x) - g'(z_y), z_x - z_y \ra
\;\; = \;\; 
\la \nabla G(x) - \nabla G(y), z_x - z_y \ra \\
\\
& = & 
\la \nabla G(x) - \nabla G(y), x - \frac{1}{\mu}\nabla G(x) - y + \frac{1}{\mu} \nabla G(y) \ra \\
\\
& = & 
- \frac{1}{\mu}\| \nabla G(x) - \nabla G(y) \|^2 + \la \nabla G(x) - \nabla G(y), x - y\ra.
\ea
$$
Applying the Cauchy-Schwartz inequality for the second term completes the proof. \qed

\subsection{Proof of Lemma~\ref{LemmaClose}}

Indeed, by the definition of smoothing, we have for any $y \in \R^n$,
$$
\ba{rcl}
G(x) & \leq & g(y) + \frac{\mu}{2}\|y - x\|^2.
\ea
$$
Substituting $y := x$, we obtain $G(x) \leq g(x)$.

On the other hand, by using that $g$ is $1$-Lipschitz, we conclude
$$
\ba{rcl}
G(x) & \Def & \min\limits_{y \in \R^n}
\Bigl\{  g(y) + \frac{\mu}{2}\|y - x\|^2   \Bigr\} \\
\\
& \geq & g(x)  +  \min\limits_{y \in \R^n}
\Bigl\{   -\|y - x\| + \frac{\mu}{2}\|y - x\|^2   \Bigr\} \\
\\
& = & g(x) - \frac{1}{2\mu}.
\ea
$$
\qed

\vspace*{-1cm}

\subsection{Proof of Lemma~\ref{LemmaLocal}}

Let us denote the objective of smoothing operator by
$$
\ba{rcl}
h(y) & \Def & g(y) + \frac{\mu}{2}\|y - x\|^2.
\ea
$$
For any point $y$ outside interior of the ball, 
we have $\|y - x\| \geq \frac{2}{\mu}$,
and hence
$$
\ba{rcl}
h(y) &  \geq & g(x) - \|y - x\| + \frac{\mu}{2}\|y - x\|^2 \\
\\
& = & g(x) + \|y - x\| \cdot \bigl( \frac{\mu}{2}\|y - x\| - 1 \bigr) \\
\\
& \geq & g(x) \;\; = \;\; h(x).
\ea
$$
So, the value of $h(\cdot)$ outside the interior is always greater than or equal
to the value at the center of the ball. Therefore, due to strong convexity, the minimum is in the interior.
\qed

\end{document}